\documentclass[11pt]{amsart}

\usepackage{parskip}
\usepackage{geometry, amsmath, amssymb, amsthm, mathrsfs, setspace, comment, amscd, latexsym, mathtools}
\usepackage[all]{xy}
\usepackage[dvipdfmx]{hyperref}

\usepackage[abs]{overpic} 
\usepackage[dvipdfmx]{pict2e}

\usepackage{color,graphicx}

\makeatletter
\@addtoreset{equation}{section}

\makeatother

\newtheorem{theorem}{Theorem}[section]

\newtheorem{corollary}[theorem]{Corollary}

\theoremstyle{definition}
\newtheorem{definition}[theorem]{Definition}

\newtheorem{remark}[theorem]{Remark}

\newcommand{\del}{\partial}

\newcommand{\Z}{\mathbb{Z}}
\newcommand{\R}{\mathbb{R}}
\newcommand{\C}{\mathbb{C}}
\newcommand{\Int}{\textrm{Int}}

\newcommand{\red}[1]{\textcolor{magenta}{\bf #1}}

\title{Surfaces in $D^4$ with the same boundary and fundamental group}
\author[Takahiro Oba]{Takahiro Oba}
\address{Department of Mathematics, Tokyo Institute of Technology, 2-12-1 Ookayama, Meguroku, Tokyo 152-8551, Japan}
\email{oba.t.ac@m.titech.ac.jp}
\thanks{This work was partially supported by JSPS KAKENHI Grant Number 15J05214.}
\date{\today}

\begin{document}

\maketitle

\begin{abstract}
We construct a family of pairs of non-isotopic symplectic surfaces in the standard symplectic $4$-disk $(D^4, \omega_{st})$ 
such that they are bounded by the same transverse knot in the standard contact $3$-sphere and 
fundamental groups of their complements are isomorphic. 
In the appendix, we prove explicitly that one can obtain a symplectic surface in $(D^4, \omega_{st})$ from 
a braided surface in a bidisk.
 \end{abstract}

\section{introduction}

This paper is concerned with 
symplectic surfaces in the standard symplectic $4$-disk $(D^4, \omega_{st})$ 
bounded by the same transverse link in the standard contact $3$-sphere $(S^3, \xi_{st})$.  
Such surfaces have been studied in some papers \cite{Ge,CGHS, Au, BV}.
Up to the present, the provided families of distinct symplectic surfaces bounded by the same transverse knot (or link) 
can be distinguished by 
the fundamental groups of their complements. 
Hence it is natural to ask whether there is a pair of non-isotopic 
symplectic surfaces in $D^4$ bounded by the same transverse knot such that 
complements of two surfaces have isomorphic fundamental groups.

The main result of this paper is the following: 

\begin{theorem}\label{theorem}
There is a family $\{(S_{1}(n), S_{2}(n))\}_{n \in \Z_{\geq0}}$ of pairs of symplectic surfaces in the standard symplectic 
$4$-disk $(D^{4},\omega_{st})$ with 
contact boundary such that:  
\begin{enumerate}
\item For a fixed $n\in \Z_{\geq 0}$, 
	\begin{enumerate}
		\item their boundaries $\del S_{j}(n)$ ($j=1,2$) are the same transverse knot up to isotopy in the boundary $(S^3, \xi_{st})$, 
		\item two fundamental groups $\pi_{1}(D^4 \setminus S_{j}(n))$ are isomorphic, and 
		\item double branched covers $X_{j}(n)$ of $D^{4}$ branched along $S_{j}(n)$ 
		are not homeomorphic, and, particularly, two surfaces $S_{j}(n)$ are not isotopic; 
	\end{enumerate}
\item The boundaries $\del S_{j}(n)$ and $\del S_{j}(n')$ are not smoothly isotopic in $\del D^4$ if $n\neq n'$.
\end{enumerate}
\end{theorem}

Rudolph exhibited two braid factorizations of a fixed $3$-braid in \cite{Ru}. 
Based on this example, we construct inequivalent braid factorizations, which provide symplectic surfaces in the above family. 
Note that it does not directly follow from the argument in \cite{Ru2} that 
from a braided surface in a bidisk 
one can obtain a symplectic surface in $(D^4, \omega_{st})$ isotopic to the given braided surface.
To complete this, we prove it in Appendix \ref{section: symplectic surfaces} explicitly.   

A \textit{Stein filling} of a contact manifold is 
a sublevel set of a proper, bounded below strictly plurisubharmonic function on a complex manifold 
whose convex boundary is contactomorphic to the given one (see \cite{Oz} for more details).
We obtain the following corollary from the above theorem combined with an argument about contact and Stein structures. 

\begin{corollary}\label{corollary}
There is a family of contact $3$-manifolds $\{ (M(n), \xi(n))\}_{n\in \Z_{\geq 0}}$ such that 
each contact manifold admits two non-homeomorphic Stein fillings $X_{1}(n)$ and $X_{2}(n)$ 
which have the same fundamental group and homology group but non-isomorphic intersection forms.
\end{corollary}

\noindent {\bf {Acknowledgements}}
The author would like to express his gratitude to Professor Hisaaki Endo for his continuous support and encouragement. 
He also thanks Dennis Auroux and Kyle Hayden for many fruitful suggestions and comments on an early draft of the paper, 
and Josh Sabloff for his kind e-mail correspondence.

\section{braided surfaces} 

\subsection{Braid groups} 

We here briefly review braid groups (see \cite[Section 2.2]{TF} for example).
Let $D^{2}$ be a closed disk in $\R^2$ equipped with the standard orientation and $K \subset \Int D^{2}$ a finite set.
Suppose that $\# K=m$. 

\begin{definition}
The \emph{braid group} with respect to $D^{2}$ and $K$, denoted by $B_{m}[D^{2},K]$, is 
the group of isotopy classes of orientation-preserving diffeomorphisms $\beta$ of $D^{2}$ such that 
$\beta|_{\del D^{2}} = id _{\del D^{2}}$ and $\beta(K)=K$. 
The elements of this group are called \emph{braids}.
\end{definition}

Let $\sigma$ be a smooth simple path in $\Int D^{2}$ with distinct end points $a,b \in K$ and $\sigma \cap K = \{ a,b\}$. 
Choose a small tubular neighborhood $U \subset \Int D^{2}$ of $\sigma$ such that 
$U \cap K = \{ a,b\}$. 

\begin{definition}
The \emph{half-twist $H(\sigma)$ along $\sigma$} is an element of the braid group $B_{m}[D^{2},K]$ 
which switches the end points $a$ and $b$ of $\sigma$ by a counterclockwise $180^{\circ}$ rotation and 
whose support is contained in $U$.
\end{definition}


\subsection{Braided surfaces and their descriptions}

Let $D_{1}^{2}$ and $D_{2}^{2}$ be two oriented closed disks. 

\begin{definition}
A \emph{braided surface} in the bidisk $D_{1}^{2} \times D_{2}^{2}$ is a properly embedded surface $S$ in $D_{1}^{2}\times D_{2}^{2}$ 
such that: 
\begin{enumerate}
\item The restriction of the first projection $pr_{1}|_{S}: S \rightarrow D_{1}^{2}$ is a simple branched covering;  
\item For each branch point $x\in S$ of $pr_{1}|_{S}$, there are complex coordinates $(z,w)$ and $\zeta$ around 
$x$ and $pr_{1}(x)$, respectively, 
compatible with orientations of $D_{1}^{2}\times D_{2}^{2}$ and $D_{1}^{2}$ such that 
$pr_{1}$ can be written as $ \zeta = pr_{1}(z,w)=z$ and locally the set $\{ (z,w)| z=w^{2}\}$ coincides with $S$.
\end{enumerate}
\end{definition}

Suppose that $S$ is a braided surface in $D_{1}^{2} \times D_{2}^{2}$. 
Let $\Delta(S) \subset \Int D_{1}^{2}$ denote the set of branch points of the covering $pr_{1}|_{S}$. 
For a point $y$ of $D_{1}^{2} \setminus \Delta(S)$, 
the number $m = \# (S \cap pr_{1}^{-1}(y))$ is called the \textit{degree} of the braided surface $S$. 

One can read off the fundamental group of the complement of a braided surface $S \subset D_{1}^{2} \times D_{2}^{2}$
from its braid monodromy. 
Fix a base point $y_{0} \in\del D_{1}^{2}$ and 
set $D_{y_{0}}= pr_{1}^{-1}(y_{0})$ and $K(y_{0})= D_{y_{0}} \cap S =\{x_{1}, \dots, x_{m} \}$.
For a point $y$ of $\Delta(S)$, 
consider a smooth simple loop $\gamma:[0,1] \rightarrow D_{1}^{2} \setminus \Delta(S)$ around $y$ 
based at $y_{0}$ 
whose bounding region does not contain any other branch points.
This loop lifts to $(\gamma([0,1]) \times D_{2}^{2})\cap S$ as a motion 
$$pr_{2}(\{ x_{1}(t), \dots, x_{m}(t) \in S\,  |\,  t \in [0,1]\})$$ 
of $m$ distinct points of $D_{2}^{2}$, 
where $pr_{2}:D_{1}^{2} \times D_{2}^{2} \rightarrow D_{2}^{2}$ is the second projection. 
When $t=0,1$, 
it is nothing but $K(y_{0})$.
Hence 
this motion defines a braid $\beta(\gamma) \in B_{m}[D_{y_{0}}, K(y_{0})]$,
called a \emph{braid monodromy} (with respect to $y_{0}$) around the branch point $y$. 
It is known that 
this braid is the half-twist $H(\sigma)$ along a smooth simple path $\sigma$ connecting two distinct points of $K(y_{0})$.
One can associate an element of $B_{m}[D_{y_{0}}, K(y_{0})]$ to any loop in $D_{1}^{2} \setminus \Delta(S)$ based at $y_{0}$, 
and define the homomorphism 
$$\varphi: \pi_{1}(D_{1}^{2} \setminus \Delta(S), y_{0}) \rightarrow B_{m}[D({y_{0}}), K(y_{0})].$$ 
Set $\Delta(S) = \{ y_{1}, \dots, y_{k}\}$. 
Take smooth simple loops $\gamma_{i} \in \pi_{1}(D_{1}^{2} \setminus \Delta(S), y_{0})$ around $y_{i}$, as we did before, 
so that 
the composition $\gamma_{1}\cdots \gamma_{k}$ is homotopic to $\del D_{1}^{2}$. 
Obviously, $\{ \gamma_{1}, \dots, \gamma_{k}\}$ serves as a free basis for $\pi_{1}(D_{1}^{2}\setminus \Delta(S), y_{0})$, 
and it is called a \emph{geometric basis} for the group.
Then, the braid $\varphi(\del D_{1}^{2})=\varphi (\gamma_{1}\cdots \gamma_{k})$ 
can be factorized into 
$k$ half-twists as 
$$\beta(\gamma_{1}) \cdots \beta(\gamma_{k}),$$ which is a \emph{braid monodromy factorization} of $\varphi(\del D_{1}^{2})$. 
As a similar notion, 
a \emph{braid factorization} of a braid $\beta$ is 
a factorization $\beta=\beta_{1} \cdots \beta_{k}$ into half-twists $\beta_{j}$. 
Remark that given a braid factorization $\beta_{1}\cdots \beta_{k}$ of a braid, 
one can construct a braided surface $S$ with $k$ branch points whose 
braid monodromy around each branch point $y_{i}$ is $\beta_{i}$ for some geometric basis for 
$\pi_{1}(D_{1}^{2} \setminus \Delta(S), y_{0})$.

\begin{figure}[t]
	\vspace{5pt}
	\centering
	\begin{overpic}[width=300pt,clip]{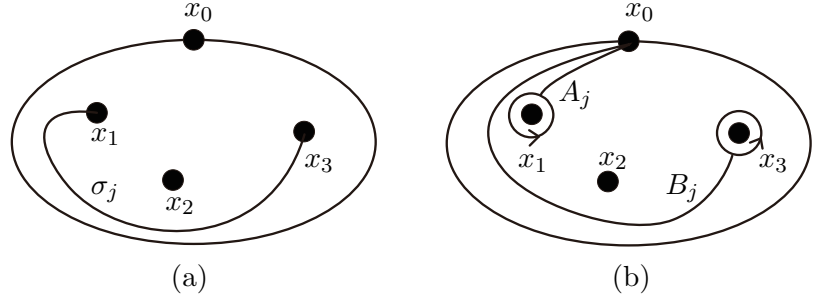}
	 \linethickness{3pt}  
  	\put(65,87){$x_{0}$}
	\put(230,87){$x_{0}$}
	\put(30,40){$x_{1}$}
	\put(190,30){$x_{1}$}
	\put(58,14){$x_{2}$}
	\put(220,31){$x_{2}$}
	\put(110,30){$x_{3}$}
	\put(280,30){$x_{3}$}
	\put(60,-15){(a)}
	\put(225,-15){(b)}
	\put(30,20){$\sigma_{j}$}
	\put(205,55){$A_{j}$}
	\put(245,20){$B_{j}$}
	\end{overpic}
	\vspace{15pt}
	\caption{(a) Path $\sigma_{j}$. (b) Loops $A_{j}$ and $B_{j}$ associated to $\sigma_{j}$.}
	\label{fig: splitting loop}
\end{figure}

Now we explain how to compute the fundamental group of the complement of a braided surface 
as the special case of \cite[Theorem 2.5]{TF}. 
Let $S$ be a braided surface in $D_{1}^{2} \times D_{2}^{2}$. 
Suppose that $\{  \gamma_{1}, \dots, \gamma_{k}\}$ is a geometric basis for the fundamental group $\pi_{1}(D_{1}^{2}\setminus \Delta(S), y_{0})$ 
and the ordered $k$-tuple   
$(H(\sigma_{1}), \dots, H(\sigma_{k}))$
consists of braid monodromies $\varphi(\gamma_{1}), \dots, \varphi(\gamma_{k})$ of $S$, 
where each $\sigma_{j}$ ($j=1, \dots, k$) is a smooth simple path connecting 
two distinct points of $K(y_{0})$. 
Fix a point $x_{0}$ of $\del D_{y_{0}}$. 
Label the points of $K(y_{0})$ as $x_{1}, \dots, x_{m}$ and 
let $\{ \gamma'_{1}, \dots, \gamma'_{m}\}$ be a geometric basis for $\pi_{1}(D_{y_{0}} \setminus K(y_{0}), x_{0})$ 
constructed in the same way we did for $\pi_{1}(D_{1}^{2} \setminus \Delta(S), y_{0})$. 
For each $j=1,\dots,k$, set 
$A_{j}=\gamma'_{i}$, where $x_{i}$ is either of end points of $\sigma_{j}$, and $B_{j}=H(\sigma_{j})(A_{j})$ (see Figure \ref{fig: splitting loop}). 
It is clear that $B_{j}$ can be expressed in terms of $\gamma'_{1}, \dots, \gamma'_{m}$ because 
they form a geometric basis for $\pi_{1}(D_{y_{0}} \setminus K(y_{0}), x_{0})$. 
By using Zariski-Van Kampen's theorem, 
we have the following formula: 
\begin{eqnarray}
\label{eqn: pi_{1}} \pi_{1}(D_{1}^{2} \times D_{2}^{2} \setminus S, x_{0}) 
 &\cong & {\pi_{1}(D_{y_{0}}\setminus K(y_{0}), x_{0})}/{ \langle\,  A_{j}=B_{j} \ (j=1,\dots, k)\,  \rangle}  \\ 
 &\cong & \langle \,  \gamma'_{1}, \dots, \gamma'_{m}\, | \,   A_{j}=B_{j} \ (j=1,\dots, k)\, \rangle.\nonumber
\end{eqnarray}
Here the point $x_{0} \in D_{y_{0}}$ is considered as one of $D_{1}^{2} \times D_{2}^{2}$ by the inclusion 
$D_{y_{0}} \hookrightarrow D_{1}^{2} \times D_{2}^{2}$.

\subsection{Double branched covers and Lefschetz fibrations}\label{section: double cover}
Let $S$ be a braided surface of degree $m$ in a bidisk $D_{1}^{2} \times D_{2}^{2}$ whose braid monodromy factorization 
with respect to some base point $y_{0}$ and geometric basis for $\pi_{1}(D_{1}^{2} \setminus \Delta(S), y_{0})$ 
is 
$$
H(\sigma_{1}) \cdots H(\sigma_{k}).  
$$
Consider the double branched covering $p:X \rightarrow D_{1}^{2} \times D_{2}^{2}$ whose branch set is $S$. 
The covering $p$ restricts to the double branched covering $p|_{F_{y_{0}}}: F_{y_{0}} = p^{-1}(D_{y_{0}}) \rightarrow D_{y_{0}}$. 
Each path $\sigma_{j}$ lifts, with respect to $p|_{F_{y_{0}}}$, 
to a unique simple closed curve $c_{j}$ on the surface $F_{y_{0}}$ up to isotopy. 
Then, according to \cite[Proposition 1]{LP}, the composition $pr_{1} \circ p: X \rightarrow D_{1}^{2}$ is 
a \emph{Lefschetz fibration} 
(see \cite[Chapter 8]{GS} for the precise definition) whose 
fibers are diffeomorphic to the surface $F_{y_{0}}$ and monodromy factorization is 
$$\tau(c_{k}) \circ \cdots \circ \tau(c_{1}).$$
Here $\tau(c)$ denotes the isotopy class of a right-handed Dehn twist along $c$. 
Throughout this paper, we use 
the functional notation for the products in the mapping class group of $F_{y_{0}}$, i.e. 
$f \circ g$ means that we apply $g$ first and then $f$.

\section{Proof of Results}

\subsection{Proof of Theorem \ref{theorem}}
\begin{figure}[t]
\vspace{5pt}
\centering
	\begin{overpic}[width=300pt,clip]{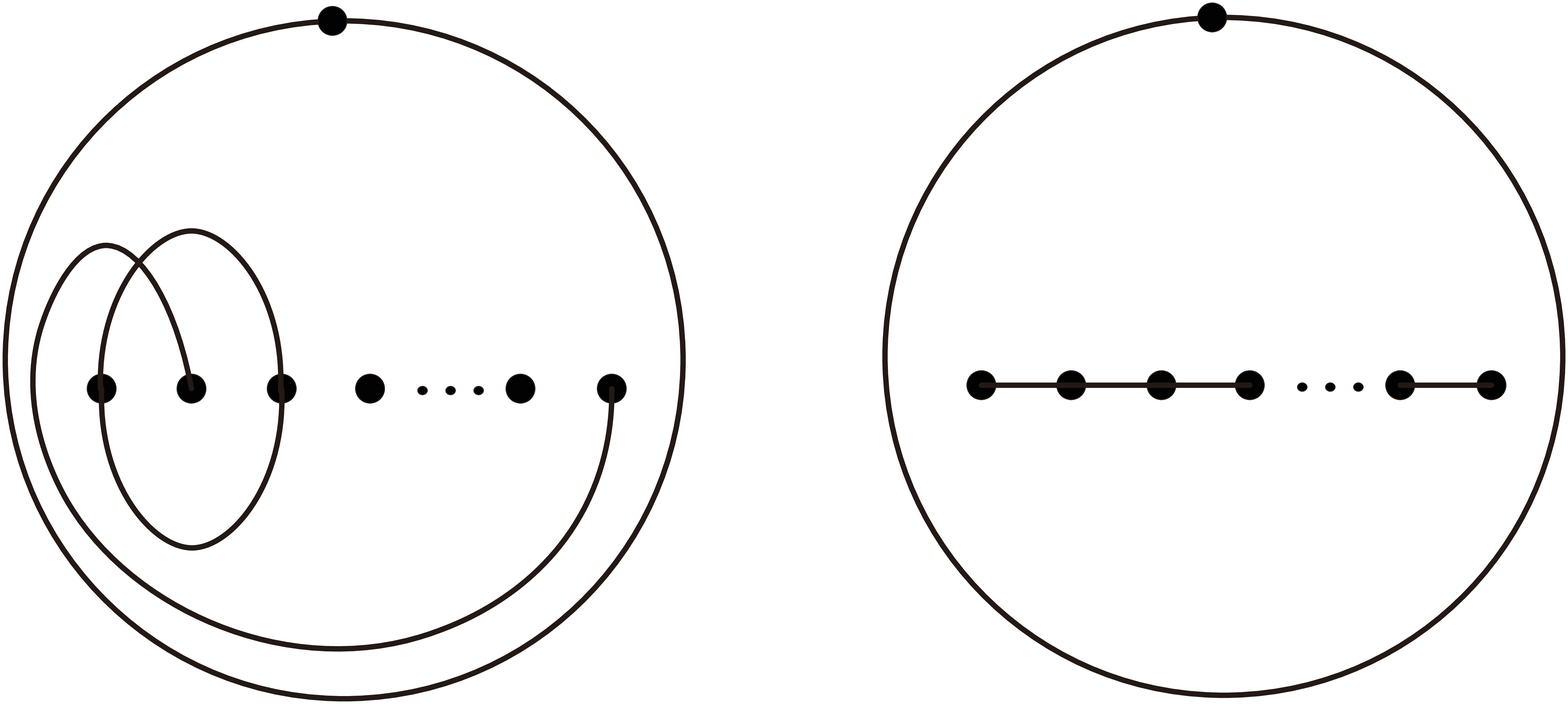}
	 \linethickness{3pt}  
  	\put(60,137){$x_{0}$}
	\put(22,50){$x_{1}$}
	\put(35,50){$x_{2}$}
	\put(55,50){$x_{3}$}
	\put(72,50){$x_{4}$}
	\put(92,50){$x_{n+2}$}
	\put(105,67){$x_{n+3}$}
	\put(50,30){$a$}
	\put(50,83){$b$}
	\put(70,15){$c_{n}$}
	\put(230,137){$x_{0}$}
	\put(185,50){$x_{1}$}
	\put(202,50){$x_{2}$}
	\put(218,50){$x_{3}$}
	\put(235,50){$x_{4}$}
	\put(273,50){$x_{n+3}$}
	\put(192,66){$d_{1}$}
	\put(210,66){$d_{2}$}
	\put(225,66){$d_{3}$}
	\put(272,66){$d_{n+2}$}
	\end{overpic}
\caption{Arcs $a,b,c_{n}, d_{1}, \dots, d_{n}$ on the disk $\mathbb{D}^{2}$ with $K_{n+3}$}	
\label{fig: arcs}
\end{figure}

Fix an integer $n \in \Z_{\geq 0}$. 
Let $\mathbb{D}^{2}$ be the closed unit disk in $\C$ and 
$K_{n+3}$ the set of $n+3$ points of $\mathrm{Int}\, \mathbb{D}^{2}$ on the real axis. 
Let $a, b, c_{n}, d_{1}, \dots, d_{n+2}$ be smooth simple paths in 
$\mathbb{D}^{2}$ as shown in Figure \ref{fig: arcs}. 
Define two braids $\beta_{1}(n), \beta_{2}(n) \in B_{n+3}[\mathbb{D}^{2},K_{n+3}]$ with factorizations given by 
\begin{align}
\label{eqn: 1} 
\beta_{1}(n) & = 
H(a) \cdot H(b) \cdot H(d_{1}) \cdot H(c_{n}) \cdot H(d_{n+2}) \cdot \, \cdots \, \cdot H(d_{3}), \\
\label{eqn: 2} \beta_{2}(n) & = 
H(H(d_{2})(a)) \cdot H(H(d_{2})(b)) \cdot H(d_{1}) \cdot H(c_{n}) \cdot H(d_{n+2}) \cdot\,  \cdots\, \cdot H(d_{3}). 
\end{align}
When $n=0$, we set $\beta_{1}(0)= H(a) \cdot H(b) \cdot H(d_{1}) \cdot H(c_{0})$ and apply the same manner to $\beta_{2}(0)$.
We obtain two braided surfaces $S_{1}(n)$ and $S_{2}(n)$ 
whose braid monodromy factorizations are the given braid factorizations (\ref{eqn: 1}) and (\ref{eqn: 2}), respectively. 
From the argument in Appendix~\ref{section: symplectic surfaces}, 
these braided surfaces can be considered as symplectic surfaces in the $4$-disk $(D^{4}, \omega_{st})$ with 
transverse knot boudaries.

The transverse knots $\del S_{1}(n)$ and $\del S_{2}(n)$ in $(S^{3}, \xi_{st})$ are 
represented by the closure of braids $\beta_{1}(n)$ and $\beta_{2}(n)$, 
respectively. 
It can be easily checked that $H(d_{2})$ commutes with the product $H(a) \cdot H(b)$, 
which proves that $\beta_{1}(n)=\beta_{2}(n)$. 
Hence two boundaries are transversely isotopic. 
Hereafter, for the sake of simplicity, set $\beta(n) =\beta_{1}(n)=\beta_{2}(n)$.

Next, we show that the fundamental groups of complements $D^{4} \setminus S_{1}(n)$ and $D^{4} \setminus S_{2}(n)$ are isomorphic. 
Fixing a base point $x_{0}$ in a fiber of the projection $pr_{1}$ of the bidisk, 
by the formula (\ref{eqn: pi_{1}}), 
$\pi_{1}(D^{4} \setminus S_{1}(n), x_{0})$ is isomorphic to the group 
generated by $\gamma_{1}, \dots, \gamma_{n+3}$ with relations 
$$
\gamma_{1}=\gamma_{2} \gamma_{3} \gamma_{2}^{-1}, \ \ \gamma_{1}=\gamma_{3}, \ \ \gamma_{1}=\gamma_{2}, 
$$
$$
(\gamma_{1} \gamma_{2} \cdots \gamma_{n+2}) \gamma_{n+3} (\gamma_{1} \gamma_{2} \cdots \gamma_{n+2})^{-1} = \gamma_{2}, \ \ 
\gamma_{j} = \gamma_{j+1} \ (j=3, \dots, n+2). 
$$
Hence $\pi_{1}(D^{4} \setminus S_{1}(n), x_{0}) \cong \langle \gamma_{1}| - \rangle \cong \Z$.
On the other hand, 
$\pi_{1}(D^{4} \setminus S_{2}(n), x_{0})$ is isomorphic to the group 
generated by $\gamma_{1}, \dots, \gamma_{n+3}$ with relations 
$$
\gamma_{2}=\gamma_{3}^{-1} \gamma_{2}^{-1} \gamma_{1} \gamma_{2} \gamma_{3}, \ \ \gamma_{1}=\gamma_{2} , \ \ \gamma_{1}=\gamma_{2}, 
$$
$$
(\gamma_{1} \gamma_{2} \cdots \gamma_{n+2}) \gamma_{n+3} (\gamma_{1} \gamma_{2} \cdots \gamma_{n+2})^{-1}=\gamma_{2}, \ \ 
\gamma_{j} = \gamma_{j+1} \ (j=3, \dots, n+2). 
$$
Thus, $\pi_{1}(D^{4} \setminus S_{2}(n), x_{0}) \cong \langle \gamma_{1}| - \rangle \cong \Z$ that 
is isomorphic to $\pi_{1}(D^{4} \setminus S_{1}(n), x_{0})$. 

\begin{figure}[t]
\vspace{5pt}
\centering
	\begin{overpic}[width=150pt,clip]{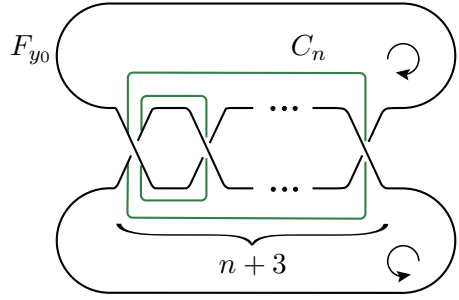}
	 \linethickness{3pt}  
	 \put(-17,90){$F_{y_{0}}$}
	 \put(88,90){$C_{n}$}
	 \put(61,8){$n+3$}
	\end{overpic}
	\caption{Surface $F_{y_{0}}$ as the double branched cover and lift $C_{n}$: Each rounded arrow indicates the orientation of $F_{y_{0}}$.}
	\label{fig: lift}
\end{figure}

\begin{figure}[t]
\vspace{15pt}
\centering
	\begin{overpic}[width=400pt,clip]{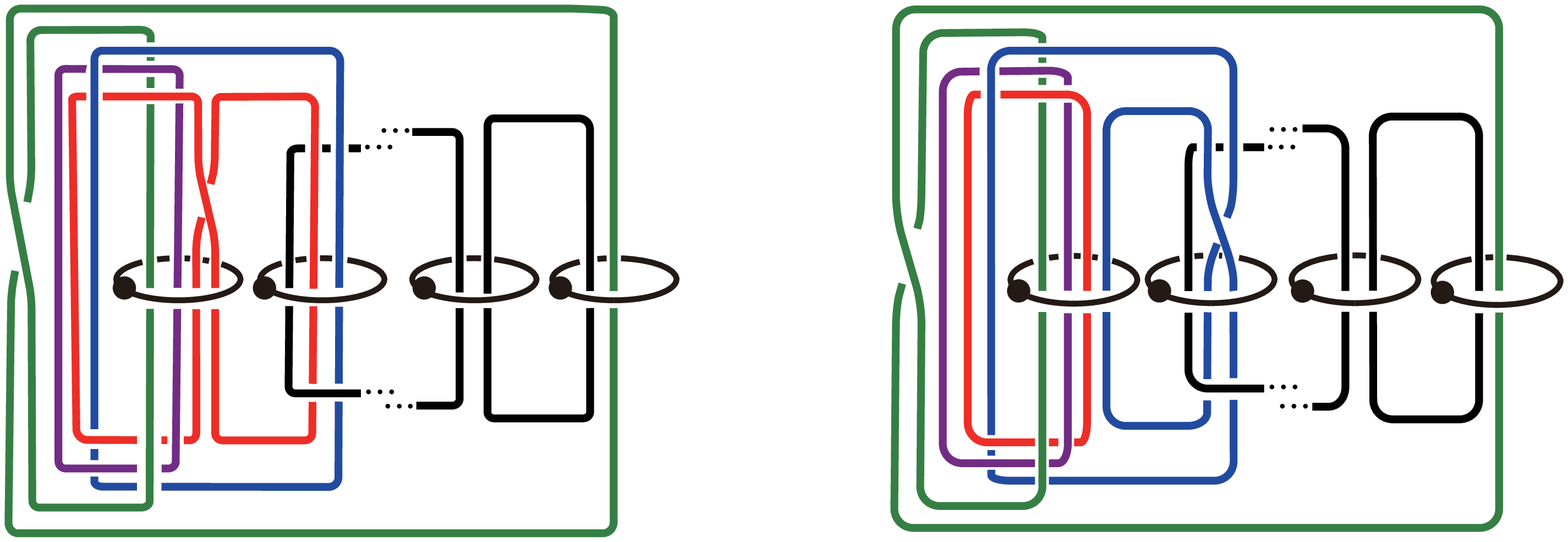}
  	\put(140,150){$-4$}
	\put(55,127){$-4$}
	\put(70,7){$X_{1}(n)$}
	\put(300,7){$X_{2}(n)$}
	\put(90,29){$C_{n}$}
	\put(60,39.5){$A$}
	\put(60,52){$B$}
	\put(13,12){$D_{1}$}
	\put(58,70){$D_{3}$}
	\put(98,132){$D_{n+1}$}	
	\put(122,136){$D_{n+2}$}
	\thicklines
	\put(20,22){\vector(0,1){18}}
	\put(365,150){$-4$}
	\put(316,145){$-4$}
	\put(362,31){$C_{n}$}
	\put(317,39){$\tau(D_{2})(A)$}
	\put(190,12){$\tau(D_{2})(B)$}
	\put(242,12){$D_{1}$}
	\put(288,70){$D_{3}$}
	\put(326,133){$D_{n+1}$}
	\put(352,137){$D_{n+2}$}
	\put(248,22){\vector(0,1){20}}
	\put(220,22){\vector(1,1){28}}
	\end{overpic}
\caption{Handle diagrams of $X_{1}(n)$ and $X_{2}(n)$: All $2$-handle framings which are not written here are $-2$.}	
\label{fig: Kirby diagram}
\end{figure}

For each $j=1,2$ let $p_{j}(n): X_{j}(n) \rightarrow D_{1}^{2} \times D_{2}^{2}$ be the double branched covering 
whose branch set is $S_{j}(n)$. 
As we discussed in Section \ref{section: double cover}, 
$X_{j}(n)$ is considered as the total space of the Lefschetz fibration $f_{j}(n)=pr_{1}\circ p_{j}(n)$. 
Let $A, B, C_{n}, D_{i}$ be lifts of arcs $a, b,c_{n}, d_{i}$, respectively, with respect to the covering $p_{j}(n)|_{F_{y_{0}}}$, 
where $F_{y_{0}}$ is the preimage of $y_{0} \in D_{1}^{2}$ under $f_{j}(n)$
(see Figure~\ref{fig: lift}). 
Fibers of the Lefschetz fibration $f_{j}(n)$ are diffeomorphic to $F_{y_{0}}$ and 
its monodromy factorization is 
\begin{align}
\label{monodromy: 1} & \tau(D_{3}) \circ \cdots \circ  \tau(D_{n+2}) \circ \tau(C_{n}) \circ \tau(D_{1}) \circ \tau(B) \circ \tau(A)  & \text{if $j=1$,} 
\\
\label{monodromy: 2} & \tau(D_{3}) \circ \cdots \circ  \tau(D_{n+2}) \circ \tau(C_{n}) \circ \tau(D_{1}) \circ \tau(\tau(D_{2})(B)) \circ 
\tau(\tau(D_{2})(A))  & \text{if $j=2$}. 
\end{align}
From these data, one can draw handle diagrams (or Kirby diagrams) of $X_{1}(n)$ and $X_{2}(n)$ as in Figure~\ref{fig: Kirby diagram}. 
Here we use the standard Seifert surface for the $(2,n+3)$-torus link as the fiber surface 
to see the monodromy curves more easily. 
We should note that the surface framing of each curve does not always coincide with its blackboard framing 
(see \cite[Section 6.3]{GS}, which explains the way to draw handle diagrams of Milnor fibers in the same manner as ours). 
After sliding $2$-handles and cancelling $1$-/$2$-handle pairs as indicated in Figure~\ref{fig: X_{1}(n)} and \ref{fig: X_{2}(n)}, 
we obtain handle diagrams of $X_{1}(n)$ and $X_{2}(n)$ each of which consists of only one $0$-handle and two $2$-handles.

From the bottom left diagram of Figure~\ref{fig: X_{2}(n)} one can see that 
$X_{2}(n)$ contains a smooth surface with self-intersection number $-2$. 
In contrast, 
we will show below that the double cover $X_{1}(n)$ contains no such surfaces: 
Let $\{ e_{1}, e_{2}\}$ be the basis for the homology group $H_{2}(X_{1}(n);\Z)$, 
where each $e_{j}$ is the homology class represented by the $2$-handle depicted in the
bottom left digram of Figure~\ref{fig: X_{1}(n)}.

\begin{figure}[t]
\vspace{10pt}
\centering
	\begin{overpic}[width=330pt,clip]{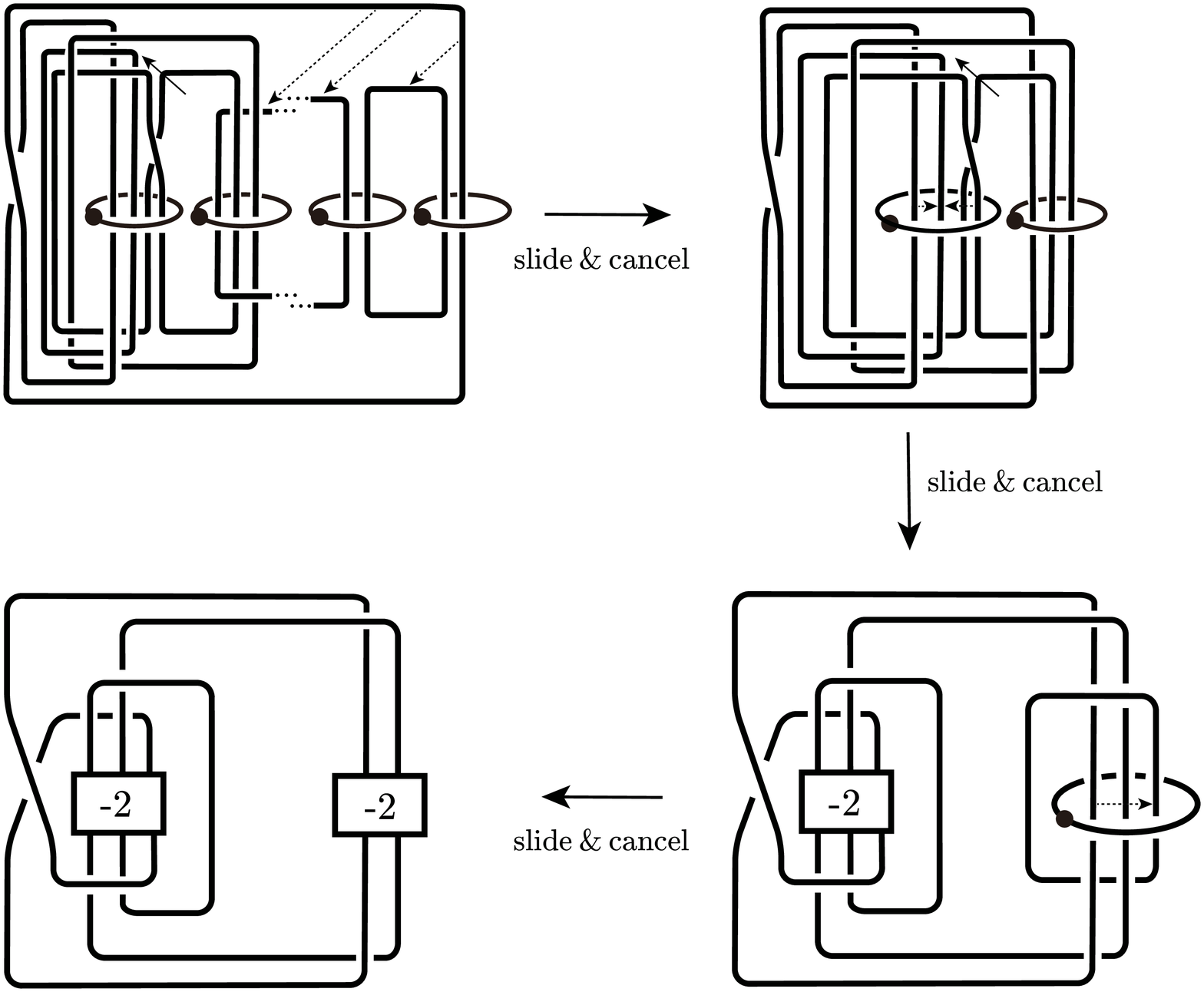}
	\put(44,240){$-2$}
	\put(45,183){$-4$}
	\put(71,258){$-2$}
	\put(72,193){$-2$}
	\put(82,178){$-2$}
	\put(102,188){$-2$}
	\put(110,164){$-4$}
	\put(269,240){$-2$}
	\put(269,183){$-4$}
	\put(286,160){$-2n-4$}
	\put(295,258){$-2$}
	\put(286,111){$-2n-4$}
	\put(280,92){$-8$}
	\put(270,20){$-2$}
	\put(80,111){$-2n-4$}
	\put(70,92){$-8$}
	\put(50,112){$e_{1}$}
	\put(50,92){$e_{2}$}
	\end{overpic}
\caption{Handle calculus for $X_{1}(n)$: Each dashed arrow in the diagram indicates how we slide a $2$-handle over another one.}	
\label{fig: X_{1}(n)}

\vspace{25pt}
	\begin{overpic}[width=330pt,clip]{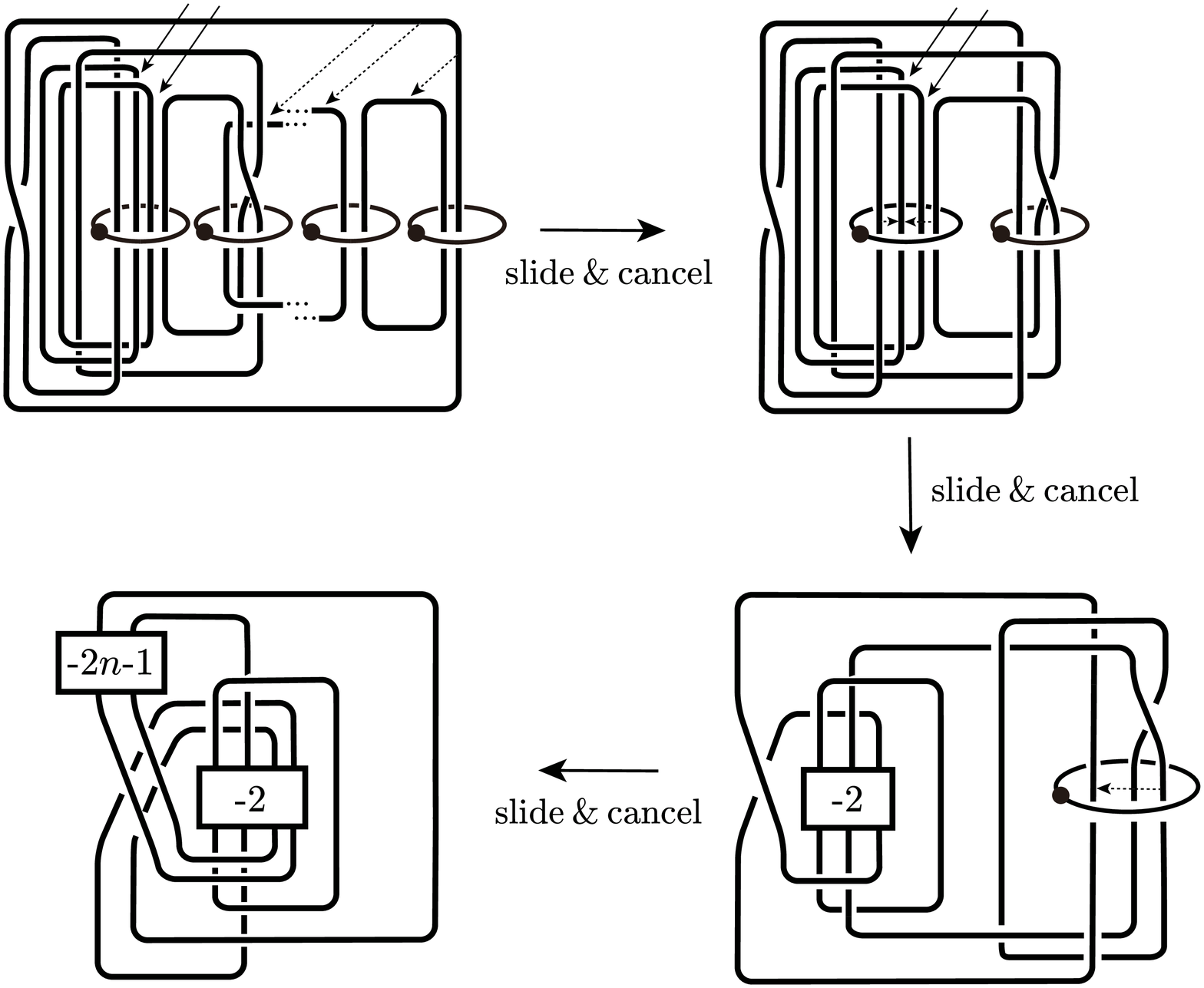}
	\put(38,262){$-2$}
	\put(53,262){$-2$}
	\put(71,183){$-2$}
	\put(77,168){$-2$}
	\put(68,160){$-4$}
	\put(100,177){$-2$}
	\put(100,155){$-4$}
	\put(242,261){$-2$}
	\put(257,261){$-2$}
	\put(282,180){$-4$}
	\put(271,150){$-2n-4$}
	\put(270,107){$-2n-4$}
	\put(273,80){$-8$}
	\put(230,22){$-2$}
	\put(80,107){$-8n-20$}
	\put(91,20){$-2$}
	\end{overpic}
\caption{Handle calculus for $X_{2}(n)$: Each dashed arrow in the diagram indicates how we slide a $2$-handle over another one.}	
\label{fig: X_{2}(n)}
\end{figure}
\clearpage

\noindent

Suppose for the sake of contradiction that there are 
integers $\alpha_{1},\alpha_{2} \in \Z$ such that 
$(\alpha_{1}e_{1}+\alpha_{2}e_{2})^{2}=-2$.
The matrix $Q(n)$ of the intersection form $Q_{X_{1}(n)}$ with respect to 
the basis $\{ e_{1}, e_{2}\}$
can be read off from the handle diagram of $X_{1}(n)$, 
and 
\begin{equation*}
Q(n)=
\begin{bmatrix} 
	-2n-4 & -1 \\ 
	-1 & -8
\end{bmatrix}.
\end{equation*}
Then, by using this matrix, 
we have 
$(-2n-4)\alpha_{1}^{2} -2 \alpha_{1} \alpha_{2} -8\alpha_{2}^{2}=-2.$
The left-hand side also can be written as the form 
$$
-(2n+4)(\alpha_{1}+\frac{1}{2n+4}\alpha_{2})^{2} - (8- \frac{1}{2n+4})\alpha_{2}^{2}. 
$$
Since the above two terms are non-positive, 
we have $- (8- {1}/{(2n+4)})\alpha_{2}^{2} \geq -2$, or $(8-{1}/{(2n+4)})\alpha_{2}^{2} \leq 2$. 
The coefficient $8- {1}/{(2n+4)}$ is greater than $2$, and hence $\alpha_{2}^2 <1$, that is, $\alpha_{2}=0$. 
Thus 
\begin{equation}
\label{eqn: 3} -(2n+4)\alpha_{1}^{2}=-2 \ \ \text{and}\ \  \alpha_{1} \in \Z \setminus \{ 0\}. 
\end{equation}
However, 
$-(2n+4)\alpha_{1}^{2} \leq -(2n+4) \leq -4$ for $\alpha_{1} \in \Z \setminus \{ 0\}$, 
which contradicts the equation (\ref{eqn: 3}). 
Thus, we conclude that $X_{1}(n)$ and $X_{2}(n)$ are not homeomorphic.

To distinguish the closure $\widehat{\beta(n)}$ from $\widehat{\beta(n')}$ for $n \neq n'$, 
we use the \emph{determinant} of a knot, defined by $|\det (V+V^{T})|$, where $V$ is a Seifert matrix for the knot. 
It is known that it equals the order of the first homology group of 
the double branched cover of $S^3$ branched along the knot. 
Moreover, let $X$ be a compact $4$-manifold admitting a handle decomposition with only one $0$-handle and 2-handles. 
Then, the determinant of a matrix for the intersection form $Q_{X}$ coincides with $|H_{1}(\del X; \Z)|$ up to sign 
(see \cite[Corollary 5.3.12]{GS}). 
Thus the determinant of the closure $\widehat{\beta(n)}$ is 
$$\det Q(n) = \det 
\begin{bmatrix} 
	-2n-4 & -1 \\ 
	-1 & -8
\end{bmatrix}
= 16n + 31,
$$
which proves that all $\widehat{\beta(n)}$ ($n \in \Z_{\geq 0}$) are mutually non-isotopic. 
This finishes the proof.

\begin{remark}
Braid factorizations $H(a) \cdot H(b)$ and $H(H(d_{2})(a)) \cdot H(H(d_{2})(b))$ we used above are essentially 
found by Rudolph in \cite[Example 1.13]{Ru}, where he showed that two factorizations are ones of the same braid. 
This pair was also used in \cite[Proposition 3.2]{Au} (see also Example 4.2 in the same paper). 
\end{remark}

\begin{remark}
Two mapping class factorizations (\ref{monodromy: 1}) and (\ref{monodromy: 2}) are related by a \textit{partial conjugation}, 
twisting the last two factors by $\tau(D_{2})$. 
This implies that two corresponding double covers are related by a Luttinger surgery along a torus built by parallel transport of 
the curve $D_{2}$ along a loop in $D_{1}^{2}$ (see \cite{ADK}).
\end{remark}

\subsection{Proof of Corollary \ref{corollary}}

Let each $S_{j}(n)$ ($j=1,2$) be the braided surface constructed above and  
$p_{j}(n): X_{j}(n) \rightarrow D^{4}$ ($j=1,2$) the double branched covering whose branch set is $S_{j}(n)$. 
As we mentioned before, the covering $p_{j}(n)$ induces the Lefschetz fibration $f_{j}(n)$ on $X_{j}(n)$. 
According to \cite{AO, LP}, $X_{j}(n)$ admits a Stein structure, 
and the contact structure $\xi_{j}(n)$ on the boundary $M_{j}(n)=\del X_{j}(n)$ 
induced from the Stein structure  
is compatible with 
the open book determined by the Lefschetz fibration. 
We see that the monodromy $\phi_{j}(n)$ of this open book is isotopic to the composition 
$$\phi_{j}(n)=
\begin{cases}
\tau(D_{3}) \circ \cdots \circ  \tau(D_{n+2}) \circ \tau(C_{n}) \circ \tau(D_{1}) \circ \tau(B) \circ \tau(A)    & \text{if $j=1$,} \\
\tau(D_{3}) \circ \cdots \circ  \tau(D_{n+2}) \circ \tau(C_{n}) \circ \tau(D_{1}) \circ \tau(\tau(D_{2})(B)) \circ \tau(\tau(D_{2})(A))   & \text{if $j=2$}. 
\end{cases}
$$
Since $\tau(D_{2})$ commutes with $\tau(B) \circ \tau(A)$, we have 
$$\tau(B) \circ \tau(A)=\tau(\tau(D_{2})(B)) \circ \tau(\tau(D_{2})(A)).$$ 
Hence $\phi_{1}(n)=\phi_{2}(n)$, which proves that the contact manifolds $(M_{1}(n),\xi_{1}(n))$ and 
$(M_{2}(n),\xi_{2}(n))$ are 
mutually contactomorphic. 
Therefore, $X_{1}(n)$ and $X_{2}(n)$ serve as Stein fillings of the contact manifold $(M(n),\xi(n)) \coloneqq (M_{1}(n),\xi_{1}(n))$ 
whose intersection forms are non-isomorphic by 
Theorem~\ref{theorem}. 
Moreover, from handle diagrams depicted in Figure~\ref{fig: X_{1}(n)} and~\ref{fig: X_{2}(n)}, 
it is obvious that $X_{1}(n)$ and $X_{2}(n)$ are simply connected and have the same homology group.
This completes the proof.

\appendix
\section{From braided surfaces to symplectic surfaces}\label{section: symplectic surfaces}
In this appendix, we explain how to obtain a symplectic surface in 
the standard symplectic $4$-disk $(D^{4}, \omega_{st})$ from a braided surface in a bidisk. 
Here $D^{4}=\{ (z,w) \in \C^{2}\, | \, |z|^{2}+|w|^{2} \leq 1 \} \subset \C^{2}$ and 
$\omega_{st}={\sqrt{-1}}(dz \wedge d\bar{z}+ dw \wedge d\bar{w})/2|_{D^4}$.
Although in fact it has been used as a well-known fact implicitly, 
here for future use we prove it explicitly. 
We assume that the reader is familiar with basics of contact and symplectic geometry. 
If necessary, we refer the reader to \cite{MS}.

Let $D^{2}(r)$ be the closed disk in $\C$ of radius $r$ centered at the origin with the complex orientation. 
In particular, $D^{2}(1)$ is the closed unit disk $\mathbb{D}^{2}$.
One can construct a braided surface $S$ in a bidisk whose braid monodromy factorization 
coincides with a given braid factorization. 
Furthermore, by a result of Rudolph \cite[Section 4]{Ru2}, for the braided surface $S$, there is 
a polynomial $f(z,w) \in \C[z,w]$ such that 
$S=\{ f(z,w)=0\} \cap (\mathbb{D}^{2} \times D^{2}(r))$ is a braided surface (with respect to the first projection of 
$\mathbb{D}^{2} \times D^{2}(r)$)
isotopic to $S$ for some $r>0$. 
Without loss of generality, we may assume that $S$ is smooth by slightly perturbing $f$. 
Moreover, 
possibly after making   
$D^{2}(r)$ larger by $D^{2}(r) \ni w \mapsto {Kw} \in D^{2}(Kr)$ 
for some large $K\gg 0$ and 
considering $\{ f(z,Kw)=0\} \cap (\mathbb{D}^{2} \times D^{2}(Kr))$, 
the boundary $\del S$ can be assumed to be a transverse link in 
$(\del \mathbb{D}^{2} \times D^{2}(r), \alpha|_{\del \mathbb{D}^{2} \times D^{2}(r)})$, 
where $\alpha= {\sqrt{-1}}(zd\bar{z}-\bar{z}dz+wd\bar{w}-\bar{w}dw)/4$.

\begin{figure}[ht]
\vspace{5pt}
	\begin{overpic}[width=140pt,clip]{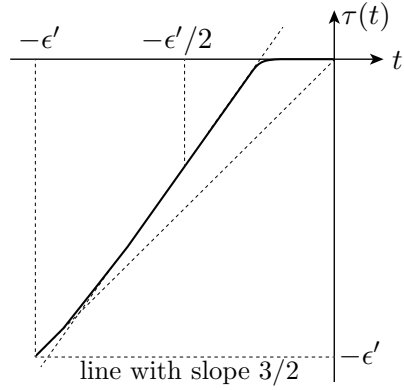}
	\put(3, 126){$-\epsilon'$}
	\put(123,7){$-\epsilon'$}
	\put(49, 126){$-\epsilon'/2$}
	\put(143, 118){$t$}
	\put(124, 135){$\tau(t)$}
	\put(25,5){\vector(-1,1){8}}
	\put(26,1){\small{line with slope $3/2$}}
	\end{overpic}
\caption{Graph of the function $\tau: (\epsilon', 0] \rightarrow \R$}	
\label{fig: graph}
\end{figure}

In addition, we may assume that the Liouville vector field 
$V=z{\del}/{\del z}-\bar{z}{\del}/{\del \bar{z}}+w{\del}/{\del w}-\bar{w}{\del}/{\del \bar{w}}$ 
for the symplectic form $\omega_{st}$
is tangent to the symplectic surface $S$.  
This can be achieved as follows: 
Since $\del S$ is a transverse link, 
there is a neighborhood $N(\del S)$ of $\del S$ in $\mathbb{D}^{2} \times D^{2}(r)$ 
endowed with a symplectomorphism 
$$\varphi: (N(\del S), \omega) \rightarrow ((-\epsilon, 0] \times S^{1} \times \R^{2}, d(e^{t}(d\theta+xdy-ydx)))$$ 
for some $\epsilon$, 
where $(t, \theta, x,y) \in (-\epsilon, 0] \times S^{1} \times \R^{2}$ and
 $\varphi(\del S) =\{ 0 \} \times S^{1} \times \{0\}$. 
The push-forward map $\varphi_{*}$ maps 
the vector field $V$ to ${\del}/{\del t}$ on $(-\epsilon, 0] \times S^{1} \times \R^{2}$. 
Then, one can find an embedding given by 
$$\psi: (-\epsilon, 0] \times S^{1} \rightarrow 
(-\epsilon, 0] \times S^{1} \times \R^{2}, \ 
\psi(t,\theta) = (t, \theta, x(t,\theta), y(t, \theta))$$ 
whose image coincides with 
$\varphi(S \cap N(\del S))$.
Note that $\psi(0,\theta)=(0, \theta, 0, 0)$ for any $\theta\in S^{1}$.
Since $\psi((-\epsilon, 0] \times S^{1})=\varphi(S \cap N(\del S))$ is symplectic, 
$$
\psi^{*}d(e^{t}(d\theta+xdy-ydx)) = (1+xy_{\theta}-x_{\theta}y+2(x_{t}y_{\theta}-x_{\theta}y_{t}))dt \wedge d\theta >0,
$$
where $x_{t}={\del x}/{\del t}$, $x_{\theta}={\del x}/{\del \theta}$, 
$y_{t}={\del y}/{\del t}$, $y_{\theta}={\del y}/{\del \theta}$. 
Now, we perturb $\psi((-\epsilon, 0] \times S^{1})$ so that ${\del}/{\del t}$ is tangent to 
the surface near the boundary. 
Since $x(t, \theta)$, $y(t, \theta)$ and their partial derivatives are continuous on $(-\epsilon,0] \times S^{1}$, 
there exists $\epsilon' \in (0, \epsilon)$ such that 
$$|x(t,\theta)y_{\theta}(t,\theta)-x_{\theta}(t,\theta)y(t,\theta)|< 1/4, \ \ 
|x_{t}(t,\theta)y_{\theta}(t,\theta) - x_{\theta}(t,\theta) y_{t}(t,\theta)|< 1/4$$
for any $(t, \theta) \in (-\epsilon',0] \times S^{1}$. 
Let $\tau: (-\epsilon', 0] \rightarrow \R$ be a smooth function 
such that (see Figure~\ref{fig: graph}): 
\begin{itemize}
\item $\tau(t)=0$ near $t=0$; 
\item $\tau(t)=t$ near $t=-\epsilon'$; 
\item $\tau'(t) \equiv 3/2$ near $t=-\epsilon'/2$; 
\item $0 \leq \tau'(t) \leq 3/2$ for any $t \in (-\epsilon', 0]$.
\end{itemize}
Using this function, 
set 
$$\psi_{s}(t, \theta)=(t, \theta,x(s\tau(t)+(1-s)t, \theta), y(s\tau(t)+(1-s)t, \theta))$$ 
for $s \in [0,1]$ as a perturbation of 
$\psi((-\epsilon', 0] \times S^{1})$. 
Then, 
$$
\psi_{s}^{*}d(e^{t}(d\theta+xdy-ydx)) = (1+xy_{\theta}-x_{\theta}y+2(s\tau'+(1-s))(x_{t}y_{\theta}-x_{\theta}y_{t}))dt \wedge d\theta,
$$
and hence all we have to do is to see that the coefficient is positive. 
Indeed, 
\begin{eqnarray*}
1+(xy_{\theta}-x_{\theta}y)+2(s\tau'+(1-s))(x_{t}y_{\theta}-x_{\theta}y_{t})  > 
1-1/4 + 2 \cdot 3/2\cdot (-1/4) =  0, 
\end{eqnarray*} 
which shows that 
$S_{s}=((\varphi^{-1} \circ {\psi_{s}})((-\epsilon', 0] \times S^{1})) \cup 
(S \setminus ((\varphi^{-1} \circ {\psi_{s}})((-\epsilon', 0] \times S^{1})))$ is symplectic.  
In particular, $S_{1}$ satisfies the desired condition because  
$\psi_{1}(t,\theta)= (t, \theta, 0,0)$ near $t=0$. 
From now on, we think of the surface $S_{1}$ as $S$.

Symplectically embed the bidisk $\mathbb{D}^{2} \times D(r)$ into the symplectic manifold $(\C^{2}, \omega)$, 
where
$\omega={\sqrt{-1}}(dz \wedge d\bar{z}+ dw \wedge d\bar{w})/2$. 
Using the flow of the Liouville vector field $V$, extend the surface $S$ and obtain the \textit{completion} $\hat{S} \subset \C^{2}$ 
of $S$, 
that is, $S \cup ([0, \infty) \times \del S)$. 
For a sufficient large $R$, 
let $S'$ be the intersection of the surface $\hat{S}$ and the round $4$-ball $D^{4}(R)$ of radius $R$. 
Then, $S'=\hat{S} \cap D^{4}(R)$ is a symplectic surface in $(D^{4}(R), \omega|_{D^{4}(R)})$. 
From this it follows that  
$$(\iota_{V}\omega)_{p}(v) = \omega_{p}(V_{p}, v)>0 $$ 
at any point $p\in \del S'$
for a tangent vector $v \in T_{p}\del S'$ such that the vector space $\langle v \rangle$ 
spanned by $v$ is isomorphic to $T_{p}\del S'$ 
as oriented vector spaces.  
This implies that the boundary $\del S'$ is a transverse link in the boundary $(\del D^{4}(R), \ker(\alpha|_{\del D^{4}(R)}))$, where 
$\alpha=\iota_{V}\omega$.

Define a diffeomorphism $\Phi_{R}: D^{4}=D^{4}(1) \rightarrow D^{4}(R)$ by $\Phi_{R}(z,w)=(Rz, Rw)$. 
Since $S'$ is symplectic in $(D^{4}(R), \omega|_{D^{4}(R)})$, 
so is $S''=\Phi_{R}^{-1}(S')$ in $(D^{4}, (\Phi_{R})^{*}\omega)$. 
The restriction of $\Phi_{R}$ on ${\del D^{4}}$ gives a contactomorphism between 
$(\del D^{4}, \ker ((\Phi_{R})^{*}\alpha|_{\del D^{4}})$ and $(\del D^{4}(R), \ker(\alpha|_{\del D^{4}(R)}))$. 
Thus, the boundary $\del S'' = \Phi_{R}^{-1}(\del S')$ is a transverse link in $(\del D^{4}, \ker ((\Phi_{R})^{*}\alpha|_{\del D^4}))$.
Now, we have 
$(\Phi_{R})^{*}\omega=R^{2}\omega$ and $(\Phi_{R})^{*}\alpha=R^{2}\alpha$. 
The coefficients $R^2>0$ do not affect whether the surface $S''$ (resp. its boundary $\del S''$) is 
symplectic (resp. transverse), and so 
$S''$ is a symplectic surface in the standard symplectic $4$-disk $(D^{4}, \omega|_{D^{4}}=\omega_{st})$ and 
$\del S''$ is a transverse link in 
the standard contact $3$-sphere
$(\del D^{4}, \ker (\alpha|_{\del D^{4}})=\xi_{st})$. 
Note that our procedure here preserves the link type of the transverse link in the boundary.
Therefore, $S'' \subset (D^{4}, \omega_{st})$ is a desired symplectic surface.

 \begin{remark}
We would like to point out that in his recent paper \cite{Hay}, 
Hayden 
obtains a similar result 
of \textit{Stein quasipositive} links in general Stein fillable contact $3$-manifolds. 
The proof is based on \textit{ascending surfaces}, originated from the work of Boileau and Orevkov \cite{BO}. 
\end{remark}

\addcontentsline{toc}{chapter}{Bibliography}
\bibliographystyle{abbrv}
\bibliography{symplecticsurfaces}

\end{document}